\documentclass[12pt,a4paper]{article}

\usepackage{amssymb,amsmath,graphicx,url,xcolor,todonotes}

\DeclareMathOperator*{\Argmax}{Arg\,max}

\outer\def\proclaim #1. #2\par{\medbreak \noindent{\bf#1.\enspace}{\sl#2}\par
	\ifdim\lastskip<\medskipamount \removelastskip\penalty55\medskip\fi}

\title{VADU 2018 Open Problem Session}

\author{Bui Thi Hoa\thanks{CIAO, Federation University Australia},  Scott B. Lindstrom\thanks{School of Mathematical and Physical Sciences, University of Newcastle, Australia}  and Vera Roshchina\thanks{School of Mathematics and Statistics, UNSW and School of Science, RMIT University, Australia} \textsuperscript{*}}

\begin{document}

\maketitle

\begin{abstract} We state the problems discussed in the open problem session at Variational Analysis Down Under (VADU2018) conference held in honour of Prof. Asen Dontchev's 70th birthday on 19--21 February 2018 at Federation University Australia, \url{https://sites.rmit.edu.au/asen/}. 
\end{abstract}
\tableofcontents

\section{Existence of local calm selections}
This problem was proposed by Asen Dontchev. All background material, including notation, history, etc. can be found in \cite{book}. We are grateful to Asen for providing this description.

\proclaim Theorem (Bartle-Graves  (1952)). Let $X$ and $Y$ be
Banach spaces and let $f:X \to Y$ be a function which is strictly
differentiable at $\bar{x}$ and such that the derivative $D f(\bar{x})$ is
surjective.  Then there exist a neighborhood $V$ of $ f(\bar{x})$ and 
a constant $\gamma > 0$ such that   $f^{-1}$ has a continuous selection $s$  on $V$ which is calm with constant $\gamma$; that is,
$$
\|s(y) - \bar{x}\| \leq \gamma\|y - f(\bar{x})\|
\;\text{for every }  y \in V.
$$

When $X$ and $Y$ are finite dimensional, even Hilbert, the proof is easy. For Banach spaces, the proof is highly nontrivial. A generalization of the Bartle-Graves theorem to set-valued mappings was obtained in \cite{a}.

Here is the open problem:

\proclaim Conjecture.
Consider  a function $f:\mathbb{R}^n \to \mathbb{R}^m$ which is Lipschitz continuous around $\bar{x}$ 
and suppose that  all matrices $A$ in Clarke's generalized Jacobian of $f$ at $\bar x$   are surjective. Then
$f^{-1}$ has a continuous
local selection around $\bar{y}$ for $\bar{x}$ which is calm at  $\bar y = f(\bar x)$.

If $n=m$ the conjecture reduces to Clarke's inverse function theorem. For $m \leq n$, according to a theorem by Pourciau \cite{P}, under the same condition the function $f$ is metrically regular. This last result was generalized recently to Banach spaces in \cite{r}.

\section{Are $6$-polytopes $3$-linked?}
This problem was presented by Bui Thi Hoa.

A graph $G$  is $k$-linked if for any selection of $k$ pairs of all distinct vertices $Y:=\{(s_1,t_1),\ldots,(s_k,t_k)\}$, $(k \ge 1)$ there exist $k$ disjoint paths, connecting the $k$ pairs of points in $Y$. If the graph of a polytope is {\it $k$-linked} we say that the polytope is also {\it $k$-linked}.

Recall that a  $d$-polytope is a $d$-dimensional polytope, i.e. the linear span of the polytope is a $d$-dimensional space. The initial question is whether or not every $d$-polytope is {\it $\lfloor d/2 \rfloor$-linked}. And the negative answer was given by Gallivan (see \cite{Gal}) with a construction of a $d$-polytope which is not {\it $\lfloor 2(d+4)/5\rfloor$-linked}.
It had been already proven that $4$-polytopes and $5$-polytopes are {\it $2$-linked} (see \cite{Tho}, \cite{Sey}), meanwhile not all $8$-polytopes are {\it $4$-linked}. The remaining question is that if all the $6$-polytopes are {\it $3$-linked}.
\section{Is FFS3 polytope decomposable?}
This problem was suggested by David Yost, and communicated during the open problem section by Scott Lindstrom and Vera Roshchina.

A polytope is called decomposable \cite{PY} if it can be represented as Minkowski sum  of dis-similar convex bodies. Two polytopes are similar if one can be obtained from the other by a dilation and a translation.

\begin{figure}[ht]
{\centering \includegraphics[width=0.3\textwidth]{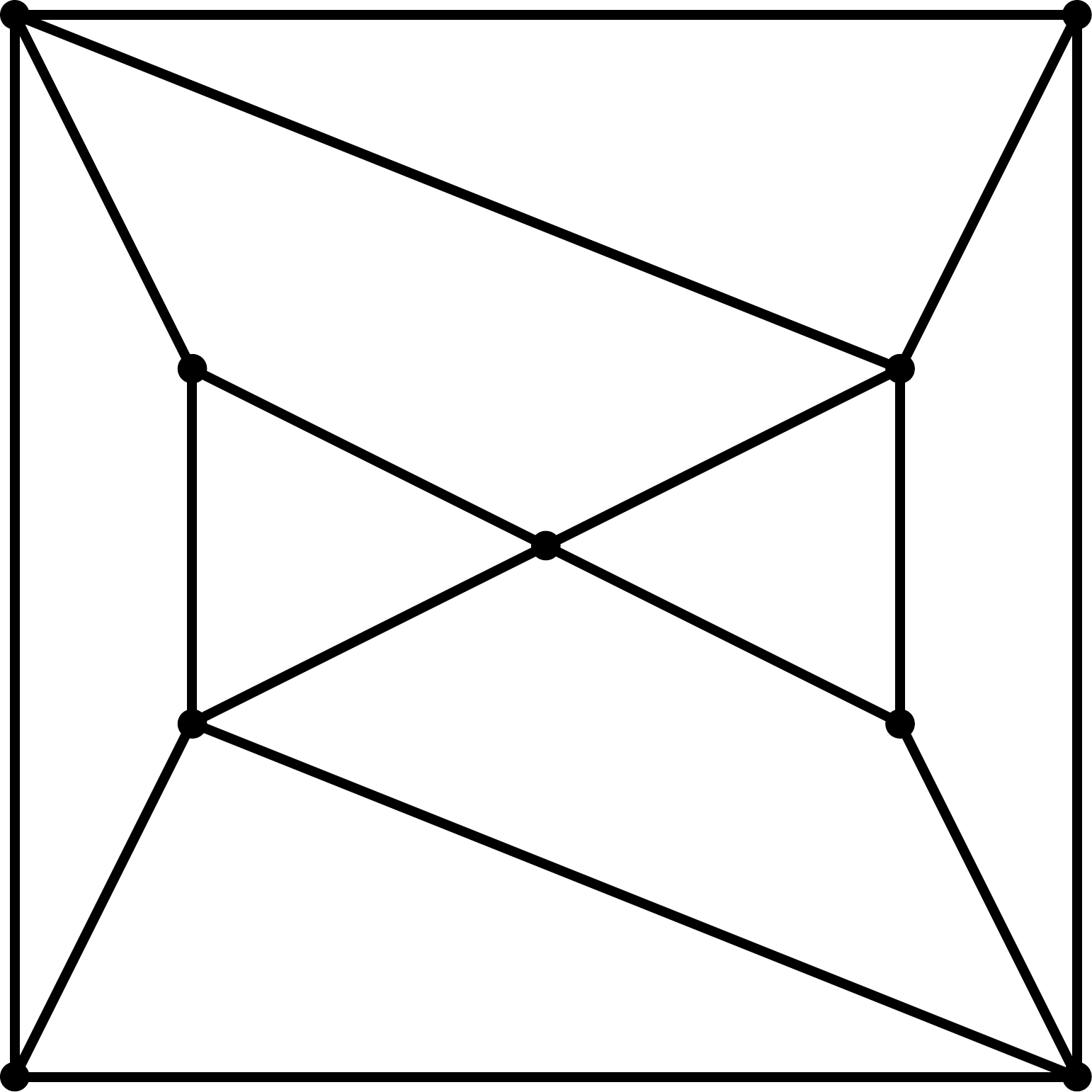}\\}
\caption{Graph of the polytope FFS3}
\label{fig:FFS3}
\end{figure}

David Yost in collaboration with Debra Briggs have classified all but one 3-polytopes with up to 16 edges in terms of decomposability (manuscript in preparation). The only remaining case is the (combinatorial) polytope FFS3 with its graph shown in Fig.~\ref{fig:FFS3}.
It is conjectured that this  polytope has no decomposable geometric realisation. 

All polytopes with up to 15 edges are classified in terms of their decomposability \cite{briggs}, and the resolution of the decomposability question for FFS3 polytope will settle the 16-edge case. However further case-by-case decomposability classification of polyhedra with higher number of edges presents a tedious challenge, and a more interesting question is developing an algorithm to check decomposability. We note that in an overwhelming number of cases indecomposability can be checked using combinatorial conditions from \cite{PY}.

\section{Projections onto compact convex sets}
This problem was proposed by Andrew Eberhard.

Let $C_1$ and $C_2$ be compact convex sets in a Hilbert space $\mathcal{H}$. The conjecture states that there always exists a point $x\in \mathcal{H}$ such that for each of its projections $p_i$ onto $C_i$, $i\in \{1,2\}$ the relevant normals $x-p_1$ and $x-p_2$ define the hyperplanes that strongly expose the faces $\{p_1\}$ and $\{p_2\}$ of $C_1$ and $C_2$ respectively. 

Recall  (see \cite[Definition 8.27]{Fabian}) that a point $x\in C$ is strongly exposed by a linear functional $f$ if $f(x) = \sup_{x'\in C} f(x')$ and $x_k\to x$  for all sequences $\{x_k\} \subset  C$ such that $\lim f(x_k) = \sup_{x\in C}f(x)$.

\section{Convergence of the continuous time Douglas-Rachford algorithm}
This problem was proposed by Scott Lindstrom.

For the feasibility problem of finding a point in the nonempty intersection $A\cap B \ne \emptyset$ of proximal sets $A$ and $B$, the Douglas-Rachford method for a given starting point $x_0$ generates a sequence
\begin{equation*}
x_n \in Tx_{n-1}:=\left(\lambda(2P_B-{\rm Id})(2P_A-{\rm Id})+(1-\lambda){\rm Id}\right)x_{n-1}
\end{equation*}
where $P_A,P_B$ denote the usual projection operators for $A,B$ respectively and $\lambda \in (0,1]$ is usually taken to be $1/2$. When $A,B$ are also convex, the sequence $(x_n)_{n\in \mathbb{N}}$ converges weakly to a fixed point of the method (see \cite{LM} and \cite{BCL}).

\begin{figure}[ht]
	\begin{center}
		\includegraphics[width=.45\textwidth]{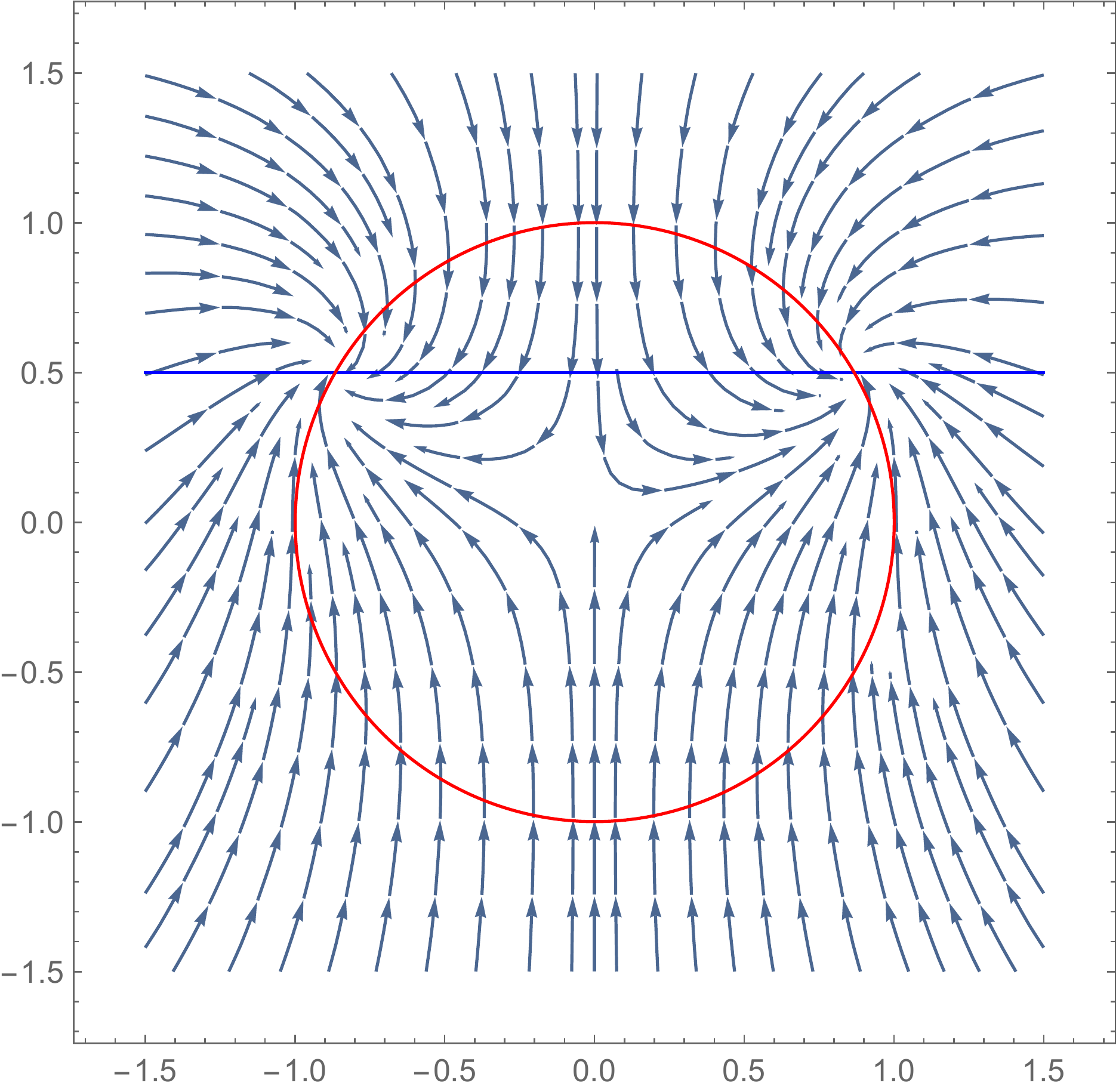}\quad \includegraphics[width=.515\textwidth]{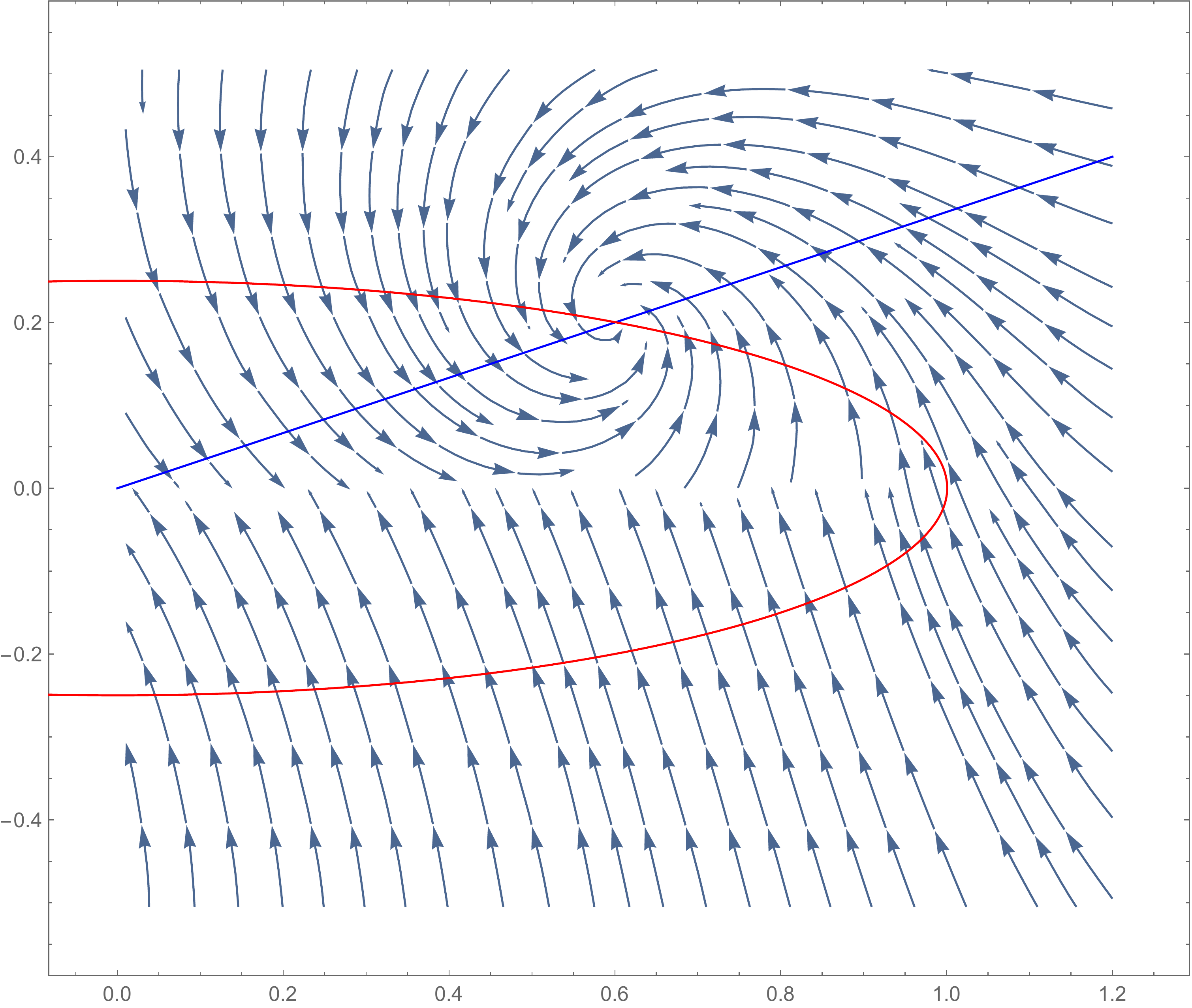}
	\end{center}
	\label{fig:circleflow}
	\caption{The flowfield \eqref{DR_DE} with a circle/line (left) and ellipse/line (right). Images courtesy of Veit Elser.}
\end{figure}

For the nonconvex case where $A$ is a circle and $B$ a line, Borwein and Sims \cite{BS} considered the ``continuous time'' version of the algorithm---whose flow field is shown at left in Figure~\ref{fig:circleflow} and corresponds to the solution of the differential equation given by
\begin{equation}\label{DR_DE}
\frac{dx}{dt}=T(x) \quad \text{when} \; \lambda \rightarrow 0^+
\end{equation} 
---as a means to approaching the question of convergence in the usual case of $\lambda=1/2$ given subtransversality, a case Benoist \cite{Benoist} answered in the affirmative by means of a Lyapunov function and which has since been extended by Minh N. Dao and Matthew K. Tam \cite{DT}. 

\begin{figure}
	\begin{center}
		\includegraphics[angle=90,width=\textwidth]{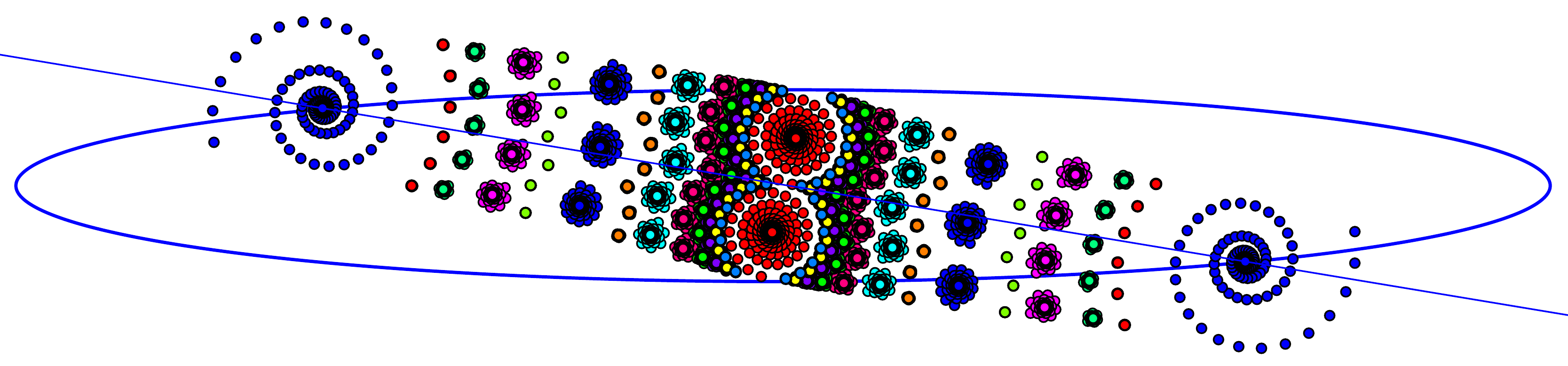}
	\end{center}
	\caption{Behaviour of Douglas-Rachford method with an ellipse and line varies from the case of a circle and line.}\label{fig:ellipseandline}
\end{figure}

The generalization to a subtransversal ellipse and line and also to a p-sphere and a line was considered by Borwein et al.\cite{BLSSS}, who showed that local convergence remains while global behaviour becomes far more complicated. See, for example, Figure~\ref{fig:ellipseandline}. Veit Elser has suggested analysing the continuous time version of the method in these more general settings and has generously furnished the images in Figure~\ref{fig:circleflow}.

\section{Minimal distance problem}
This problem was proposed by Alex Kruger.

Given a finite set of points $a_1,\dots, a_m\in X$, where $X$ is an Euclidean space, find the solution to the problem
\begin{equation}\label{eq:pb-alex}
\min_{x\in X}\max_{i\in \{1,\dots, m\}}\|a_i-x\|.
\end{equation} 
The problem has a unique solution for which $x$ is the centre of the minimal Euclidean sphere that contains all points. However it is unclear whether there exists a neat way to write this explicitly.

This is a particular case of a more general problem. The space $X$ can be an arbitrary normed linear or even a metric space. In the latter case, the norm of the difference in \eqref{eq:pb-alex} should be replaced by the distance. Instead of the maximum in \eqref{eq:pb-alex}, it could be an arbitrary norm in $\mathbb{R}^m$.

\section{Demyanov-Ryabova conjecture}
This problem was communicated by Vera Roshchina.

The problem was originally stated in \cite[Conjecture 1]{DR}. Recently two different special cases were confirmed in \cite{DP,TS}. During the preparation of this file a counterexample was found \cite{Vera-DR}.

Given a  finite family $\Omega$ of convex polytopes in $\mathbb{R}^n$, for each unit vector $g\in S_{n-1}$ we construct a new polytope as the convex hull of all support faces of all polytopes in the family $\Omega$, i.e. we define the function 
$$
C(g) := \mathrm{conv\,} \{\Argmax_{x\in P}\langle x, g\rangle \,|\,  P \in \Omega\}.
$$
Collecting all such polytopes, we obtain a new finite family of polytopes,
$$
F(\Omega) = \{C(g)\, g\in S_{n-1}\}.
$$

Now starting from a given finite collection of polytopes $\Omega_0$ we apply this transformation infinitely obtaining a sequence $\Omega_0$, $\Omega_1$, $\Omega_2$, \dots, where $\Omega_i = F(\Omega_{i-1})$, $i\in \mathbb{N}$.

The original Demyanov-Ryabova conjecture claimed that this sequence eventually reaches a two-cycle, i.e. for a sufficiently large $N$ we have $\Omega_{N+2} = \Omega_N$. Since we now know that the conjecture is false, the question is to find a characterisation of such collections of polytopes that yield two-cycles, extending and generalising the results of \cite{DP,TS}.

\section{D\"urer's conjecture}
This problem was communicated by Vera Roshchina.

\begin{figure}[ht]
{\centering
\includegraphics[scale=1]{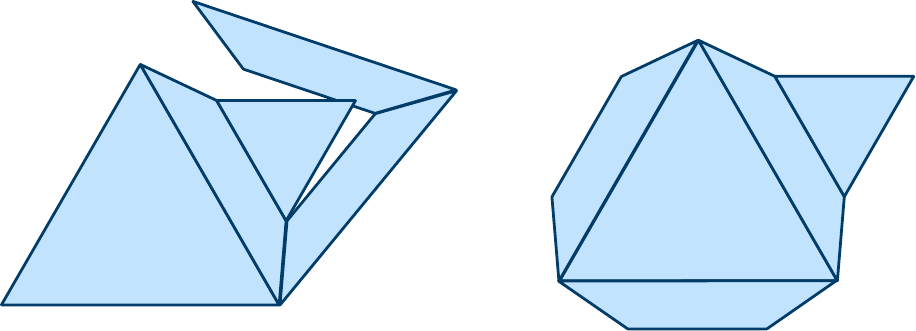}\\}
\caption{Two different nets of the same truncated tetrahedron}
\label{fig:unfolding}
\end{figure}

Albrecht D\"urer dedicated a nontrivial part of his career to laying out the geometric foundations of drawing and perspective. His five centuries old work \cite{Durer} is available online via Google books. The mathematical statement known as D\"urer's conjecture was motivated by this work and proposed in 1975 by Shephard \cite{Shephard}. A net (or unfolding) of a 3-polytope is the process of cutting it along its edges,  so that the resulting connected shape can be flattened (developed) into the plane \cite{GhomiNotices}. It is not difficult to find examples of polytopes for which certain cuts result in overlapping unfoldings, such as the truncated tetrahedron shown in Fig.~\ref{fig:unfolding}
(see \cite{GhomiProof}). 

The D\"urer's conjecture is a claim that any polytope has a nonoverlapping net. A significant recent contribution in this direction is the work by Mohammed Ghomi who showed that every polytope is combinatorially equivalent to an unfoldable one \cite{GhomiProof}. For more details we refer the reader to an overview \cite{GhomiNotices} by the same author.

\section*{Acknowledgements}

We are grateful to Asen Dontchev, Andrew Eberhard, Alex Kruger and David Yost for patiently clarifying the mathematical details of their open problems to us. 

\bibliographystyle{plain}
\bibliography{references}

\begin{thebibliography}{10}

\bibitem{BCL}
H.H. Bauschke, P.L. Combettes, and R.D. Luke.
\newblock Phase retrieval, error reduction algorithm, and fienup variants: a
  view from convex optimization.
\newblock {\em Journal of the Optical Society of America A}, 19:1334--1345,
  2002.

\bibitem{Benoist}
Joel Benoist.
\newblock The douglas-rachford algorithm for the case of the sphere and line.
\newblock {\em Journal of Global Optimization}, 63:363--380, 2015.

\bibitem{BLSSS}
Jonathan~M. Borwein, Scott~B. Lindstrom, Brailey Sims, Matthew Skerritt, and
  Anna Schneider.
\newblock Dynamics of the douglas-rachford method for ellipses and p-spheres.
\newblock {\em Set Valued and Variational Analysis}.

\bibitem{BS}
Jonathan~M. Borwein and Brailey Sims.
\newblock The douglas-rachford algorithm in the absence of convexity.
\newblock In Heinz~H. Bauschke, Regina~S. Burachik, Patrick~L. Combettes, Veit
  Elser, D.~Russell Luke, and Henry Wolkowicz, editors, {\em Fixed Point
  Algorithms for Inverse Problems in Science and Engineering}, volume~49.
  Springer Optimization and Its Applications, 2011.

\bibitem{briggs}
Debra Briggs.
\newblock Comprehensive catalogue of polyhedra.
\newblock {\em AMSI Report}, 2016.

\bibitem{Tho}
Thomassen C.
\newblock 2--linked graphs.
\newblock {\em European J. Combin.}, 1:371--378, 1980.

\bibitem{r}
Radek Cibulka, Asen~L. Dontchev, and Vladimir~M. Veliov.
\newblock Lyusternik-graves theorems for the sum of a lipschitz function and a
  set-valued mapping.
\newblock {\em SIAM J. Control Optim.}, 54(6):3273--3296, 2016.

\bibitem{Sey}
Seymour~P. D.
\newblock Disjoint paths in graphs.
\newblock {\em Discrete Math.}, 29:293--309, 1980.

\bibitem{DP}
Aris Daniilidis and Colin Petitjean.
\newblock A partial answer to the demyanov-ryabova conjecture.
\newblock {\em Set-Valued and Variational Analysis}, Jul 2017.

\bibitem{DT}
Minh~N. Dao and Matthew~.K. Tam.
\newblock A lyapunov-type approach to convergence of the douglas-rachford
  algorithm.
\newblock 2017.

\bibitem{DR}
Vladimir~F. Demyanov and Julia~A. Ryabova.
\newblock Exhausters, coexhausters and converters in nonsmooth analysis.
\newblock {\em Discrete Contin. Dyn. Syst.}, 31(4):1273--1292, 2011.

\bibitem{a}
Asen~L. Dontchev.
\newblock A local selection theorem for metrically regular mappings.
\newblock {\em J. Convex Anal.}, 11(1):81--94, 2004.

\bibitem{book}
Asen~L. Dontchev and Terry~R. Rockafellar.
\newblock {\em Implicit Functions and Solution mappings. A View From
  Variational Analysis}.
\newblock Springer, 2014.

\bibitem{Durer}
Albrecht D\"urer.
\newblock {\em Underweysung der Messung, mit dem Zirckel und Richtscheyt, in
  Linien, Ebenen unnd gantzen corporen}.
\newblock Nuremberg, 1525.

\bibitem{Fabian}
Mari\'an Fabian, Petr Habala, Petr H\'ajek, Vicente Montesinos Santaluc\'\i~a,
  Jan Pelant, and V\'aclav Zizler.
\newblock {\em Functional analysis and infinite-dimensional geometry}, volume~8
  of {\em CMS Books in Mathematics/Ouvrages de Math\'ematiques de la SMC}.
\newblock Springer-Verlag, New York, 2001.

\bibitem{GhomiProof}
Mohammad Ghomi.
\newblock Affine unfoldings of convex polyhedra.
\newblock {\em Geom. Topol.}, 18(5):3055--3090, 2014.

\bibitem{GhomiNotices}
Mohammad Ghomi.
\newblock D\"urer's unfolding problem for convex polyhedra.
\newblock {\em Notices Amer. Math. Soc.}, 65(1):25--27, 2018.

\bibitem{LM}
P.L. Lions and B.~Mercier.
\newblock Splitting algorithms for the sum of two nonlinear operators.
\newblock {\em SIAM Journal on Numerical Analysis}, 16:964--979, 1979.

\bibitem{P}
Brouce~H. Pourciau.
\newblock Analysis and optimization of lipschitz continuous mappings.
\newblock {\em J. Opt. Theory Appl.}, 22:311--351, 1977.

\bibitem{PY}
Krzysztof Przes{\l}awski and David Yost.
\newblock More indecomposable polyhedra.
\newblock {\em Extracta Math.}, 31(2):169--188, 2016.

\bibitem{Vera-DR}
Vera Roshchina.
\newblock Demyanov-ryabova conjecture is false, 2018.

\bibitem{Gal}
Gallivan S.
\newblock Disjoint edge paths between given vertices of a convex polytope.
\newblock {\em J. Combin.Theory Ser. A}, 39:112–--115, 1985.

\bibitem{TS}
Tian Sang.
\newblock On the conjecture by {D}emyanov-{R}yabova in converting finite
  exhausters.
\newblock {\em J. Optim. Theory Appl.}, 174(3):712--727, 2017.

\bibitem{Shephard}
G.~C. Shephard.
\newblock Convex polytopes with convex nets.
\newblock {\em Math. Proc. Cambridge Philos. Soc.}, 78(3):389--403, 1975.

\end{thebibliography}

\end{document}